\documentclass[11pt,a4paper,leqno,twoside, headinclude]{amsart}

\usepackage[tracking=true]{microtype}
\DeclareMicrotypeSet*[tracking]{my}%
  { font = */*/*/sc/* }%
\SetTracking{ encoding = *, shape = sc }{ 45 }

\title{A note on the Hilali conjecture}
\def\titl{A note on the Hilali conjecture}
\def\auth{Manuel Amann}
\date{January 12th, 2015}


\author{\auth}

\usepackage[english]{babel}

\usepackage[colorlinks,pdftex, plainpages=false]{hyperref}
\hypersetup{pdftitle=\titl, pdfauthor=\auth, pdftoolbar=false,
plainpages=false, hyperindex=true, pdfdisplaydoctitle=true}

\hypersetup{colorlinks,%
citecolor=black,%
filecolor=black,%
linkcolor=black,%
urlcolor=black,%
pdftex}

\usepackage{color}

\usepackage{amsmath, amssymb, amscd, amsthm}
\usepackage{stmaryrd}
\usepackage{latexsym}
\usepackage[all]{xy}
\usepackage{pb-diagram}
\usepackage{rotating}
\usepackage{multicol}
\usepackage{lscape}
\usepackage{wasysym}
\usepackage{longtable}
\usepackage{enumerate}

\xyoption{all}

\newtheorem{theo}{Theorem}[section]
\newtheorem{main}{Theorem}

\newtheorem*{main*}{Theorem}

\newtheorem*{mainprop*}{Proposition}

\newtheorem{mainconj}{Conjecture}

\newtheorem{defi2}[theo]{Definition}
\newtheorem*{defi2*}{Definition}

\newenvironment{defi*}{\begin{defi2*}\normalfont}{\end{defi2*}}

\newenvironment{defin*}[1]{\begin{defi2*}[#1]\normalfont}{\end{defi2*}}
\newtheorem*{rem2*}{Remark}
\newenvironment{rem*}{\begin{rem2*}\normalfont}{\hfill$\boxbox$\end{rem2*}}
\newtheorem{rem2}[theo]{Remark}
\newenvironment{rem}{\begin{rem2}\normalfont}{\hfill$\boxbox$\end{rem2}}

\newtheorem*{cor*}{Corollary}

\newtheorem*{conj*}{Conjecture}
\newtheorem*{theo*}{Theorem}
\newtheorem*{ques*}{Question}
\newtheorem*{mi2}{Main Idea}

\newtheorem{ex2}[theo]{Example}

\newtheorem{exer2}[theo]{Exercise}

\newtheorem{alg2}[theo]{Algorithm}

\newcommand{\nn}{{\mathbb{N}}}                                     
\newcommand{\qq}{{\mathbb{Q}}}                                     
\newcommand{\s}{{\mathbb{S}}}                                      
\newcommand{\dif} {{\operatorname{d}}}                             
\newcommand{\In} {{\,\subseteq\,}}                                 
\newcommand{\im} {{\operatorname{im\,}}}                           
\newcommand{\coker}{{\operatorname{coker\,}}}                        
\newcommand{\rk}{{\operatorname{rk\,}}}                            

\newcommand{\comment}[1]{}                                         
\newcommand{\xto}[1]{\xrightarrow{#1}}                             
\newcommand{\hto}[1]{\overset{#1}{\hookrightarrow}}                

\newcommand{\odd}{\textrm{odd}}                                    
\newcommand{\even}{\textrm{even}}                                  

\begin{document}

\maketitle \thispagestyle{empty}


\begin{abstract}
In this short note we observe that the Hilali conjecture holds for $2$-stage spaces, i.e.~we argue that the dimension of the rational cohomology is at least as large as the dimension of the rational homotopy groups for these spaces. We also prove the Hilali conjecture for a class of spaces which puts it into the context of fibrations.
\end{abstract}


\section*{Introduction}

There are a lot of prominent conjectures in Rational Homotopy Theory centering around the question of how large the cohomology algebra of a space has to be under certain conditions. The toral rank conjecture poses this question in the context of (almost) free torus actions, the Halperin conjecture---interpreted in the right way---asks this for the cohomology of the total space of certain fibrations. A conjecture by Hilali poses such a question based upon the size of the rational homotopy groups of a space (see \cite{Hil80}). It is formulated for rationally elliptic spaces, i.e.~for simply connected spaces $X$ satisfying $\dim \pi_*(X)\otimes \qq<\infty$ and $\dim H^*(X;\qq)< \infty$.

\begin{conj*}[Hilali]
Let $X$ be a simply-connected rationally elliptic space. Then it holds that
\begin{align*}
\dim H^*(X;\qq)\geq \dim \pi_*(X)\otimes \qq
\end{align*}
\end{conj*}

We remark that in terms of minimal Sullivan models $(\Lambda W,\dif)$ this conjecture can be stated equivalently as
\begin{align*}
\dim H(\Lambda W,\dif) \geq \dim W
\end{align*}
In this article we shall essentially use this transition.

The Hilali conjecture was confirmed in several cases. It is known if $X$ is a nilmanifold, if $\dim X\leq 16$, if $X$ is a formal or a coformal space, or if it is a symplectic or cosymplectic manifold (see \cite{NY11}, \cite{HM13}, \cite{HM08a}, \cite{HM08}, \cite{HM13a}), respectively for  hyperelliptic spaces (see \cite{FFM14}).

\vspace{5mm}

In this note we consider the two-stage case as well as a class of spaces which involves fibration constructions.

This is the class $(\mathcal I)$ of elliptic spaces, characterised by the subsequent properties of the minimal Sullivan model $(\Lambda W,\dif)$. For this we draw upon the \emph{n-stage decomposition} of $(\Lambda W,\dif)=\Lambda (\bigoplus_{0\leq i} W_i, \dif)$ with $\dif(W_i)\In \Lambda W_{< i}$ for all $i$. Besides, recall the definition of the \emph{Wang derivation} $\theta_{w_i}$ by
\begin{align*}
w_i\cdot \theta=\dif-\bar \dif_{w_i}
\end{align*}
where we define $\bar \dif_{w_i}=\dif|_{w_i=0}$ as the usual differential followed by evaluation of $w_i=0$.

Let $(\Lambda V,\dif')$ be an elliptic minimal Sullivan algebra which satisfies the Hilali conjecture and admits a $n_2$-stage decomposition for some $n_2\in \nn_0$. Let $W_{n_1+1+i}=V_i$ for $0\leq i\leq n_2$.
Let $W_0\oplus \ldots \oplus W_{n_1}\oplus W_{n_1+n_2+1} \oplus \ldots \oplus W_{n_3}$ be concentrated in odd (ordinary) degree. The class $\mathcal{I}$
is defined via the existence of a homogeneous basis $(w_i)_i$---homogeneous with respect to ordinary degree $\deg$ and lower degree $\operatorname{lowdeg}$---for $W$ satisfying
that
\begin{align*}
\dif|_{\Lambda W_{n_1,\ldots, n_2}}=\dif'
\end{align*}
and one of the following two conditions.
\begin{align}\tag{$\mathcal{I}$}
\begin{minipage}[h]{10cm}
\begin{enumerate}
\item
\begin{itemize}
\item[-] $(\dif w_i)^2=0$ for $\operatorname{lowdeg} w_i>n_1+n_2$ and
\item[-] $\theta_{w_i}^2(w_j)=0$ for $\operatorname{lowdeg} w_i\leq n_1$ and $\operatorname{lowdeg} w_j\leq n_1+n_2$,  or
\end{itemize}
\item
\begin{itemize}
\item[-] $(\dif w_i)^2|_{\Lambda (W_{>n_1})}=0$ for $\operatorname{lowdeg} w_i>n_1+n_2$ and
\item[-] $\theta_{w_i}^2=0$ for $\operatorname{lowdeg} w_i\leq n_1$.
\end{itemize}
\end{enumerate}
\end{minipage}
\end{align}
where $|_{\Lambda (W_{>n_1})}$ denotes the projection $w_i\mapsto 0$ for $i\leq n_1$.

The second class of spaces we consider are the so-called \emph{$2$-stage algebras}, which are characterised by
\begin{align}\tag{$\mathcal{S}_2$}
\begin{minipage}[h]{10cm}
\begin{enumerate}
\item $\dif(U)=0$
\item $\dif(V)\In \Lambda U $
\end{enumerate}
\end{minipage}
\end{align}
(with a priori no further restrictions on the parity of the degree.) Note that neither all two-stage spaces are hyperelliptic nor that all hyperelliptic spaces are two-stage (see Section \ref{sec03} of the Appendix). However, both classes contain pure spaces  and the proofs in both cases are similar.

\begin{main}\label{theoA}
The Hilali conjecture holds for the classes $(\mathcal{I})$ and $(\mathcal{S}_2)$.
\end{main}

\begin{cor*}
Let $(\Lambda W,\dif)$ be a Sullivan algebra with $W$ concentrated in odd degree. Suppose $W$ admits a basis $w_1, \ldots, w_n$ satisfying that
$\dif w_i$ is a monomial in the $w_j$ 
for all $1\leq i\leq n$. Then $(\Lambda W,\dif)$ satisfies the Hilali conjecture.
\end{cor*}

\begin{rem*}
We note that \cite[Estimate (9), p.~32]{FFM14} also proves the toral rank conjecture for hyperelliptic spaces. In fact, on rationally elliptic spaces $X$ it holds that $\rk_0(X)\leq - \chi_\pi(X)=\dim W^\odd-\dim W^\even$ where $(\Lambda W,\dif)$ is a minimal Sullivan model for $X$ (see \cite[Theorem 7.13, p.~279]{FOT08}). In the hyperelliptic case the estimate says that $\dim H(\Lambda W,\dif)\geq 2^{-\chi_\pi(X)}$.
\end{rem*}

\begin{rem*}
We observe that \cite[Estimate (9), p.~32]{FFM14} can also be obtained by the following elementary arguments. In the terminology of \cite[Section 4.2, 32]{FFM14} one may compare the dimension of the $E_2$-term of the Leray--Serre spectral sequence associated to the rational fibration
\begin{align*}
(\Lambda V,\dif)\hto{}(\Lambda W,\dif)=(\Lambda (V \oplus  \langle \bar y_1\rangle ),\dif)\to (\Lambda \langle \bar y_1\rangle,0)
\end{align*}
with $\dif \bar y_1=x_1\in V$, namely $2\dim H(\Lambda V,\dif)$, to the dimension of its infinity term, which is $\dim H(\Lambda W,\dif)$. Proceeding inductively directly yields the inequality.
\end{rem*}


\section{The First Class of Spaces}\label{sec01}

As usual we write $W_{\leq n_1+n_2}$ for $W_0\oplus W_1 \oplus \ldots \oplus W_{n_1+n_2}$ etc.

The main characteristic of the class $\mathcal{I}$ is indeed that it allows a decomposition into rational fibrations. In the first case we obtain the fibrations
\begin{align*}
(\Lambda W_{\leq n_1+n_2},\dif)\hto{} (\Lambda W,\dif) \to (\Lambda W_{>n_2},\bar\dif)
\intertext{and}
(\Lambda W_{\leq n_1},\dif)\hto{} (\Lambda W_{\leq n_1+n_2},\dif) \to (\Lambda W_{n_1+1,\ldots, n_1+n_2},\bar\dif)
\end{align*}
and in the second one the fibrations
\begin{align*}
(\Lambda W_{\leq n_1},\dif)\hto{} (\Lambda W,\dif) \to (\Lambda W_{>n_1},\bar\dif)
\intertext{and}
(\Lambda W_{n_1+1, \ldots,n_1+ n_2},\dif)\hto{} (\Lambda W_{\geq n_1},\dif) \to (\Lambda W_{>n_2},\bar\dif)
\end{align*}

That is, in any case the space results from two fibrations involving the space $(\Lambda V,\dif')$. In the first case this space will first play the role of the fibre in a fibration, before the resulting totla space will serve as a base space in a second fibration.

In the second case $(\Lambda V,\dif')$ will play the role of the base in the first fibration, before the new total space will act as a fibre in a second fibration.

\begin{proof}[\textsc{Proof of Theorem \ref{theoA}--Second assertion}]
Let $(\Lambda W,\dif)$ be an algebra of type $\mathcal{I}$. We shall only consider the first scenario depicted above, the second one following from completely analogous arguments.

We consider the fibration
\begin{align*}
(\Lambda W_{\leq n_1},\dif)\hto{} (\Lambda W_{\leq n_1+n_2},\dif) \to (\Lambda W_{n_1+1,\ldots, n_1+n_2},\bar\dif)
\end{align*}
and observe that we may replace it by a sequence of fibrations as follows. Let $w_1,\ldots, w_k$ be the basis elements from above spanning $W_{\leq n_1}$ in increasing lower degree. We obtain a fibration
\begin{align*}
(\Lambda \langle w_1\rangle,\dif) \hto{} (\Lambda W_{\leq n_1+n_2},\dif) \to (\Lambda (W_{n_1+1,\ldots,n_1+ n_2}\oplus \langle w_2,\ldots, w_k\rangle),\bar\dif)
\end{align*}
and consider its associated Wang sequence. From \cite[Remark 2.2, p.~195]{JL04} we recall that the connecting homomorphism of the sequence corresponds to the Wang derivation from above. More precisely, we have the following long exact sequence.
\begin{align*}
\ldots & \xto{} H^i(\Lambda W_{\leq n_1+n_2},\dif) \to H^i(\Lambda (W_{n_1+1,\ldots,n_1+ n_2}\oplus \langle w_2,\ldots, w_k\rangle),\bar\dif) \\&\xto{\theta_{w_1}^*} H^{i-\deg w_1+1}(\Lambda (W_{n_1+1,\ldots, n_1+n_2}\oplus \langle w_2,\ldots, w_k\rangle),\bar\dif)\\&\to H^{i+1}(\Lambda W_{\leq n_1+n_2},\dif)\to\ldots
\end{align*}
where it is easy to see that $\theta$ induces $\theta_{w_1}^*$ on cohomology (see \cite[p.~194]{JL04}). By assumption also the morphism $\theta^*$ satisfies $(\theta_{w_1}^*)^2=0$. We have that
\begin{align*}
H(\Lambda (W_{n_1+1,\ldots, n_1+n_2}\oplus \langle w_2,\ldots, w_k\rangle),\bar\dif)\cong\ker \theta_{w_1}^*\oplus \im \theta_{w_1}^*
\end{align*}
The Wang sequence splits as
\begin{align}\label{eqn01}
0\to(\coker \theta_{w_1}^*)^{i-\deg w_1}\to H^i(\Lambda W_{\leq n_1+n_2},\dif) \to (\ker \theta_{w_1}^*)^i\to 0
\end{align}

Consequently, every element in $\ker \theta_{w_1}^*$ contributes to the cohomology of
$H(\Lambda W_{\leq n_1+n_2},\dif)$.
Choose a complement to $\ker \theta_{w_1}^*$ isomorphic to $\im \theta_{w_1}^*$. Since $(\theta_{w_1}^*)^2=0$ every such element also contributes to the cohomology via the kernel, namely in its isomorphic form from $\im \theta_{w_1}^*$, on which $\theta_{w_1}^*$ vanishes.

This shows that
\begin{align*}
&H(\Lambda W_{\leq n_1+n_2},\dif)
\\=&\dim H(\Lambda (W_{n_1+1,\ldots, n_1+n_2}\oplus \langle w_1,\ldots, w_k\rangle),\dif) \\\geq& \dim H(\Lambda (W_{n_1,\ldots, n_1+n_2}\oplus \langle w_2,\ldots, w_k\rangle),\bar\dif)
\end{align*}
We need to argue that strict inequality holds. For this we concentrate on the contribution by the cokernel in \eqref{eqn01}. Indeed, strict equality holds unless $\coker \theta_{w_1}^*=0$ and $\theta_{w_1}^*$ is surjective. Since $\deg w_1\neq 0$, the morphism $\theta_{w_1}^*$ strictly reduces the degree. As we are dealing with rationally elliptic spaces, which have finite-dimensional cohomology, in particular, it follows that surjectivity can never hold. This proves that
\begin{align*}
&\dim H(\Lambda (W_{n_1,\ldots, n_1+n_2}\oplus \langle w_1,\ldots, w_k\rangle),\dif) \\\geq& \dim H(\Lambda (W_{n_1,\ldots, n_1+n_2}\oplus \langle w_2,\ldots, w_k\rangle),\bar\dif)+1
\end{align*}

Using the analogous arguments we proceed by considering the morphisms $\theta_{w_2}^*$, $\theta_{w_3}^*$, \ldots, $\theta_{w_k}^*$. Eventually, this yields
\begin{align}\label{eqn02}
\dim H(\Lambda (W_{\leq n_1+n_2},\dif) \geq \dim H(\Lambda W_{n_1+1,\ldots, n_1+n_2},\bar\dif)+k
\end{align}

\vspace{5mm}

Let us now deal with the fibration
\begin{align*}
(\Lambda W_{\leq n_1+n_2},\dif)\hto{} (\Lambda W,\dif) \to (\Lambda W_{>n_2},\bar\dif)
\end{align*}
This time let $w_1,\ldots, w_k$ denote a basis of $W_{>n_2}$ ordered by increasing lower degree.
Again, we shall decompose this fibration into a sequence of fibrations, this time with spherical fibres, beginning with
\begin{align*}
(\Lambda (W_{\leq n_1+n_2}\oplus \langle w_1, \ldots , w_{k-1}\rangle),\dif)\hto{} (\Lambda W,\dif) \to (\Lambda \langle w_k\rangle,\bar\dif)
\end{align*}
We consider the associated Gysin sequence in order to prove that
\begin{align}\label{eqn03}
\nonumber& \dim H(\Lambda W,\dif)\\
\nonumber=&\dim H(\Lambda W_{\leq n_1+n_2}\oplus \langle w_1, \ldots , w_{k}\rangle,\dif)\\
\nonumber \geq & \dim H(\Lambda W_{\leq n_1+n_2}\oplus \langle w_1, \ldots , w_{k-1}\rangle,\dif)+1
\nonumber \\ & \vdots
\\
\nonumber \geq & \dim H(\Lambda W_{\leq n_1+n_2},\dif)+k
\end{align}
For this we argue in an analogous way to the line of reasoning with the Wang sequence above. We only present the arguments for the first spherical fibrations, all the other ones being completely analogous.

It is essential to observe that the Euler class of the sphere bundle is given as the class $[\dif w_1]$ of $\dif w_1$ in the cohomology of the base space $(\Lambda (W_{\leq n_1+n_2}\oplus \langle w_1, \ldots , w_{k-1}\rangle),\dif)$---see \cite[Example 4, p.~202]{FHT01}. By assumption $(\dif w_1)^2=0$. The connecting homomorphism of the Gysin sequence is the cup product with the Euler class. Thus we may mimic exactly the arguments from above by considering the short exact sequence
\begin{align}\label{eqn01}
0\to(\coker ([\dif w_1]\cup (\cdot)) )^{i}\to H^i(\Lambda W,\dif) \to (\ker ([\dif w_1]\cup (\cdot)))^{i-\deg w_k}\to 0
\end{align}
This yields Inequality \eqref{eqn03}.

\vspace{5mm}

Now we combine Inequalities \eqref{eqn03} and \eqref{eqn02} together with the assumption from the assertion that
\begin{align*}
&\dim H(\Lambda W_{n_1+1,\ldots, n_1+n_2},\bar\dif)
\\=&\dim H(\Lambda V,\dif)
\\\geq & \dim \pi_*(\Lambda V,\dif)
\\=&\dim V
\\=&\dim W_{n_1+1,\ldots, n_1+n_2}
\end{align*}
in order to prove the assertion, namely
\begin{align*}
&\dim H(\Lambda W,\dif)\\ \geq &  \dim H(\Lambda W_{\leq n_1+n_2},\dif)+\dim W_{> n_2}
\\ \geq & \dim H(\Lambda W_{n_1+1,\ldots, n_1+n_2},\bar\dif)+\dim W_{\leq n_1}+\dim W_{>n_2}
\\\geq &\dim W_{n_1+1,\ldots, n_1+n_2} +\dim W_{\leq n_1}+\dim W_{>n_2}
\\= & \dim W
\end{align*}

\end{proof}

\begin{rem}
The Gysin and Wang sequences are not restricted to odd-dimensional fibres or bases. In the case of a fibration with fibre $\s^2$, the Halperin conjecture, which is confirmed in this case, however, yields that the cohomology modules of the total space splits as a product and the Hilali conjecture for it holds true trivially once the base space satisfies it. We leave it to the reader to consider variations of the presented arguments whenever the base space of the fibration is an even sphere.
\end{rem}
\begin{rem}
We only used that the \emph{cohomology classes} of the differentials and the Wang derivations vanish; thus the class $(\mathcal{I})$ may be presented in slightly larger generality.
\end{rem}

\vspace{5mm}

The proof of the corollary to Theorem \ref{theoA} is obvious now: Once $W$ is concentrated in odd degrees, the square of every monomial in the $w_i$ vanishes.

\section{The Second Class of Spaces}\label{sec02}

Let $(\Lambda (U\oplus V),\dif)$ be a $2$-stage algebra.

We essentially use the inequality
\begin{align}\label{eqn1}
\dim H(\Lambda (U\oplus V),\dif) \geq 2^{\dim V-\dim U^\even}
\end{align}
from \cite[Theorem 2.3, p.~195]{JL04}.

\begin{proof}[\textsc{Proof of Theorem \ref{theoA}--Second assertion}]
We choose a basis $(x_i)_{1\leq i\leq n}$ of $U$ and  $(y_i)_{1\leq i\leq n_1+r}$ of $V$ such that
\begin{itemize}
\item[-] $n=n_1+n_2$
\item[-] $\deg x_i$ is even for $1\leq i\leq n_1$
\item[-] $\deg x_i$ is odd for $n_1+1\leq i\leq n_1+ n_2$
\item[-] $\deg y_i$ is odd for $1\leq i\leq n_1+r$
\end{itemize}

Let us justify these statements: Passing to the associated pure model of $(\Lambda (U\oplus V),\dif)$ we make the following classical observations:
\begin{itemize}
\item[-] Since this pure model has finite dimensional cohomology if and only so has the $2$-stage model, $V$ is concentrated in odd degrees.
\item[-] Since the pure model has finite dimensional cohomology, there need to exist at least $n_1$ many basis elements of $V$ (mapping to a regular sequence in $\Lambda U$ under $\dif_\sigma$, the associated pure differential).
\end{itemize}

Thus we need to prove that
\begin{align}\label{eqn3}
\dim H(\Lambda (U\oplus V),\dif)\geq n+ n_1+r
\end{align}
By definition we have that $\dim H(\Lambda (U\oplus V),\dif)\geq n+1$ (since $\dif U=0$ and taking into account zeroth cohomology).

Every element in $\Lambda^2 U$ represents a cohomology class. Every relation imposed on this vector space comes from an element of $V$ for reasons of word-length. In other words: The cohomology vector space represented by elements from $\Lambda^2 U$ has dimension
\begin{align*}
&\dim \Lambda^2 U - \dim \im \dif|_V \cap \Lambda^2 U \\
\geq & \bigg({{n} \choose 2} + n_1 \bigg)  - (n_1+r)\\
= & \frac{ n^2 -n }{2}-r
\end{align*}
Thus the cohomology represented by elements of word-length at most two is bounded from below by
\begin{align}\label{eqn2}
1+n+\frac{ n^2 -n }{2}-r=1-r+\frac{  n^2+n}{2}
\end{align}

Let us prove that except for a few cases of $(n_1,n_2,r)$ one of the following two assertions holds
\begin{align*}
& n+ n_1+r \leq
\begin{cases}
1-r+\frac{ n^2+n}{2}  & \textrm{or}\\
2^r
\end{cases}
\\ \iff &
\begin{cases}
n_1 +2 r\leq 1+\frac{ n^2-n}{2}  & \textrm{or}\\
n+ n_1+r \leq 2^r
\end{cases}
\end{align*}
By Observation \eqref{eqn3} this will prove the result in nearly all the cases due to Inequality \eqref{eqn2} in the first case respectively due to Inequality \eqref{eqn1} (and $r=\dim V-\dim U^\even$) in the second one.

So assume that $r\geq \frac{n^2-3n}{4}$ and note that $\frac{n^2-3n}{4}\leq \frac{1}{2} + \frac{n^2-n}{4}-\frac{n_1}{2}$. Under this assumption we verify that the second condition is satisfied for $n\geq 6$. Thus we assume $n_1\leq n\leq 5$. In this case $n+n_1+r\leq 2^r$ is satisfied if $r\geq 4$. Thus we have to assume that $n_1\leq n\leq 5$ and $r\leq 3$. A computer based check then reveals that none of the conditions is satisfied only if
\begin{align*}
(n_1,n_2,r)\in \{ & (1,0,1), (1,1,1), (1,1,2), (1,2,2), (2,0,1), (2,0,2), (2,1,2),\\& (2,2,3), (3,0,1), (3,0,2), (3,0,3), (3,1,3), (4,0,2), (4,0,3)\}
\end{align*}
One may do the respective simple checks in each case, or one may cite \cite{FFM14}, which proves the Hilali conjecture in the hyperelliptic case. Indeed, we may neglect all those triples with $n_2=0$, as this is the pure case. Moreover, we may also ignore the cases with $n_2=1$, since in this case the space is necessarily hyperelliptic. Thus it remains to deal with $(n_1,n_2,r)\in \{(1,2,2), (2,2,3)\}$.

\vspace{5mm}

Let us deal with $(1,2,2)$ first. In this case $\dim \pi_*(\Lambda (U\oplus V),\dif)=\dim (U\oplus V)=6$.
Note that $\dim V^\odd=2$ and $\dim (\Lambda^2 V)^\odd=2$. Since $V$ is concentrated in odd degrees, it follows that $\dim H^\odd(\Lambda (U\oplus V),\dif)\geq 4$. Since the space has vanishing Euler characteristic, it follows that $\dim H^\even(\Lambda (U\oplus V),\dif)\geq 4$ and $\dim H(\Lambda (U\oplus V),\dif)\geq 8$.

In the case of the triple $(2,2,3)$ we argue similarly: The odd-degree subspace of the cohomology algebra represented by word-length one elements is $2$-dimensional, the odd-degree subspace represented by word-length $2$ is $4$-dimensional. Vanishing Euler characteristic implies that $\dim H(\Lambda (U\oplus V),\dif)\geq 2\cdot 6=12$.

\end{proof}

\begin{rem}
The arguments concerning the special cases in the last part of the proof generalise: Since these elements have odd degree, the elements $x_ix_j$ with $1\leq i \leq n_1$, $n_1+1\leq j\leq n$, represent linearly independent elements in cohomology. So the odd-degree cohomology is at least $n_1\cdot n_2$-dimensional and due to the vanishing Euler characteristic for $n>n_1$ the cohomology will be at least $2n_1 n_2$-dimensional.
\end{rem}


\begin{appendix} \label{sec03}
\section{Separating Examples}

We briefly sketch some algebras by which one may tell apart the different classes of spaces we consider. Denote by $(\mathcal{H})$ the class of \emph{hyperelliptic algebras} $(\Lambda W,\dif)$, i.e.~elliptic algebras satisfying that
\begin{enumerate}
\item $U= W^\even$, $V=W^\odd$
\item $\dif(U)=0$
\item $\dif(V)\In \Lambda^{>0} U \otimes \Lambda V$
\end{enumerate}

We provide three minimal Sullivan algebras having the following properties:
\begin{itemize}
\item $(\Lambda X,\dif)$ is in $(\mathcal{I})$, but neither in $(\mathcal{H})$ nor in $(\mathcal{S}_2)$.
\item $(\Lambda Y,\dif)$ is in $(\mathcal{S}_2)$, but neither in $(\mathcal{H})$ nor in $(\mathcal{I})$.
\item $(\Lambda Z,\dif)$ is in $(\mathcal{H})$, but not in $(\mathcal{S}_2)$.
\end{itemize}

We define the first algebra by $X=\langle x_1,x_2,x_3,x_4\rangle$ with $\deg x_1=\deg x_2=3$, $\deg x_3=5$, $\deg x_4=7$, $\dif x_1=\dif x_2=0$, $\dif x_3=x_1x_2$, $\dif x_4=x_1x_3$. Then, by construction, this algebra is $3$-stage (and not $2$-stage), and, since $X^\even=0$, but not all differentials are trivial, it is not hyperelliptic. A direct check shows that it is in $(\mathcal{I})$, however.

The second algebra is given by $Y=\langle y_1,y_2,y_3,y_4,y_5\rangle$ with $\deg y_1=\deg y_2=\deg y_3=\deg y_4=3$, $\deg y_5=5$, $\dif y_1=\dif y_2=\dif y_3=\dif y_4=0$, $\dif y_5=y_1y_2+y_1y_3+y_1y_4+y_2y_3+y_2y_4+y_3y_4$. Obviously, this algebra is not hyperelliptic. A direct check then also shows that it does not satisfy any of the vanishing conditions for the differential or the Wang derivation required for $(\mathcal{I})$. (For example, $(\dif y_5)^2=y_1y_2y_3y_4$.)

Define the third algebra by $Z=\langle z_1,z_2,z_3,z_4, z_5,z_6, z_7\rangle$ with $\deg z_1=\deg z_2=\deg z_3=2$, $\deg z_4=\deg z_5=\deg z_6=3$, $\deg z_7=7$,  $\dif z_1=\dif z_2=\dif z_3=\dif z_6=0$, $\dif z_4=z_1z_2$, $\dif z_5=z_1z_3$, $\dif z_4=(z_3z_4-z_2z_5)z_6$. It is $3$-stage and hyperelliptic.
\end{appendix}


\vfill

\begin{center}
\noindent
\begin{minipage}{\linewidth}
\small \noindent \textsc
{Manuel Amann} \\
\textsc{Fakult\"at f\"ur Mathematik}\\
\textsc{Institut f\"ur Algebra und Geometrie}\\
\textsc{Karlsruher Institut f\"ur Technologie}\\
\textsc{Kaiserstra\ss e 89--93}\\
\textsc{76133 Karlsruhe}\\
\textsc{Germany}\\
[1ex]
\textsf{manuel.amann@kit.edu}\\
\textsf{http://topology.math.kit.edu/$21\_54$.php}
\end{minipage}
\end{center}

\end{document}